\documentclass[intlimits]{amsart}

\providecommand{\abs}[1]{\lvert#1\rvert}

\newcommand{\defn}[1]{{\em #1}}

\newcommand{\upperh}[1]{\mathfrak{h}^{#1}}
\newcommand{\bessel}[3][]{K_{#2}#1\left({#3}\right)}
\newcommand{\cospi}[1]{\cos{\left(2 \pi {#1}\right)}}

\newcommand{\whittaker}[3][]{W^{#1}\!\left(#2,#3\right)}
\newcommand{\slz}[1]{SL(#1,\mathbb{Z})}
\newcommand{\gl}[2][, \mathbb{R}]{GL(#2#1)}
\newcommand{\iwasawa}[1]{\gl{#1} / \left(O(#1,\mathbb{R}) \cdot \mathbb{R}^\times\right)}

\newcommand{\mellinM}[1]{#1^*}

\newcommand{\gammaprod}[2]{
\def\testa{#2}\def\testb{}
\,\Gamma\!\left[
\ifx\testa\testb 
#1
\else
\genfrac{}{}{0pt}{}{#1}{#2}
\fi
\right]\!
}

\newcommand{\pochh}[2][n]{\left(#2\right)_{#1}}

\newcommand{\eqn}[1]{(\ref{#1})}
\newcommand{\gammaf}[1]{\,\Gamma\!\left(\frac{#1}{2}\right)\!}

\newcommand{\iinfint}[3]{\int_{-\infty}^\infty \! \! \int_{-\infty}^\infty #1 \, d#2 \, d#3}
\newcommand{\oneint}[2]{\int_0^{1} #1 \, d#2}

\newcommand{\infint}[3][0]{\int_{#1}^{\infty} #2 \, d#3}
\newcommand{\mellin}[3]{\int_0^{\infty} #1 #2^{#3} \,\frac{d#2}{#2}}

\newcommand{\invmellin}[3][(\sigma)]{\frac{1}{2 \pi i}\int_{#1} #2 \, d#3}
\newcommand{\invmmellin}[4][(\sigma)(\tau)]{\frac{1}{(2 \pi i)^2}\iint_{#1} #2 \, d#3 \, d#4}
\newcommand{\conj}[1]{\overline{#1}}

\usepackage{latexsym}
\usepackage{amsthm}
\usepackage{url} 
\usepackage{enumerate} 

\usepackage[pdftex,plainpages=false,pdftitle={Evaluating Whittaker functions and Maass forms for SL(3,Z)},pdfauthor={Borislav Mezhericher},pdfpagemode={UseOutlines},letterpaper,bookmarks,bookmarksopen=true,pdfstartview={FitH},bookmarksnumbered=true]{hyperref}

\title{Evaluating Whittaker functions and Maass forms for $\slz{3}$}
\author{Borislav Mezhericher}
\begin{document}
\begin{abstract}
We present and compare several algorithms for evaluating Jacquet's Whittaker functions for $\slz{3}$. The most suitable algorithm is then applied to the problem of evaluating a Maass form for $\slz{3}$ with known eigenvalues and Fourier coefficients.
\end{abstract}
\maketitle
\section{Introduction}
Our goal is to provide tools for numerical experimentation with Maass forms for $\slz{3}$. The hope is that the information obtained from these experiments will be useful in formulating and testing hypotheses, which later can be rigorously proved. Even in the absence of such rigorous proofs, the data obtained by experimentation can lead to useful insights and intuition.

The main computational tool is an explicit version of the Fourier expansion of a Maass form for $\slz{3}$,  given by Goldfeld in~\cite{goldfeld}. This expansion, derived in greater generality by Piatetski-Shapiro in~\cite{piat-shapiro} and independently by Shalika in~\cite{shalika}, utilizes generalized Whittaker functions introduced by Jacquet in~\cite{jacquet}. Evaluating these functions efficiently presents a major difficulty in evaluating Maass forms. Surprisingly, with the exception of the work of K. Broughan \cite{broughan} where he computes Jacquet's Whittaker functions for $\gl{n}$ using Stade's integral representation from~\cite{stade}, we are not aware of any computational results in that area. We address that problem in Section~\ref{whittaker}, where we present several algorithms for evaluating Jacquet's Whittaker function for $\slz{3}$. One algorithm is determined to be particularly suited for evaluating a Maass form using its Fourier expansion.

A possible reason for the scarcity of results related to evaluating higher-rank Whittaker functions is that, until recently, there were no explicitly known examples of generic Maass forms for $\slz{3}$. This has changed last year:  first examples of $L$-functions of degree 3 were discovered by C. Bian and A. Booker \cite{booker-ams}, \cite{bian} followed by hundreds of other examples discovered by D. Farmer, S. Koutsoliotas and S. Lemurell \cite{farmer-etal}. According to the so-called converse theorem, the Dirichlet coefficients of each of these $L$-functions should correspond to the Fourier coefficients of a Maass for $\slz{3}$.

With these recent developments, an algorithm for evaluating a Maass form given the coefficients would be very useful for further experiments. One such experiment is given in Section~\ref{maass}, where we evaluate the Maass form that corresponds to one of the $L$-functions found by Bian and provide informal evidence of its automorphy.

\section{Preliminaries}
We need to recall some facts concerning Barnes-type integrals, automorphic forms for $\slz{3}$, and numerical integration.

\subsection{Barnes-type integrals}
Most functions that are of interest to us can be represented by a \defn{Barnes-type integral}, that is, an inverse Mellin transform of ratios of gamma functions:
\begin{equation}
\label{barnes-integral}
f(z) = \invmellin{\gammaprod{b_1+s,\dots,b_m+s,1-a_1-s,\dots,1-a_n-s}{1-b_{m+1}-s,\dots,1-b_q-s,a_{n+1}+s,\dots,a_p+s} z^{-s}}{s},
\end{equation}
where we adopt the notation
$$
\gammaprod{a_1,a_2,\ldots,a_n}{b_1,b_2,\ldots,b_k} = \frac{\Gamma(a_1)\Gamma(a_2)\cdots\Gamma(a_n)}{\Gamma(b_1)\Gamma(b_2)\cdots\Gamma(b_k)}.
$$
The path of integration is taken to be a vertical line from $(\sigma - i \infty)$ to $(\sigma + i \infty)$, possibly indented so that all poles of $\Gamma(1-a_j+s)$, $j=1,2,\dots,n$ lie to the left of the path, and all poles of $\Gamma(b_i-s)$, $i=1,2,\dots,m$ lie to the right of the path. For a description of other possible contours of integration and corresponding conditions on $a_j$'s and $b_i$'s for the integral to converge, see~\cite{luke} (where such an integral is called a \defn{G-function}).

One special case is the \defn{K-Bessel function}, defined below (compare to equation (11) in section 6.4 of~\cite{luke}):
\begin{equation}
\label{mellin-bessel}
4 \bessel{\mu}{2 \pi y} = \invmellin{\gammaprod{\frac{s+\mu}{2},\frac{s- \mu}{2}}{} (\pi y)^{-s}}{s}.
\end{equation}
By moving the line of integration to the right, we can see that the K-Bessel function decays rapidly as $y\rightarrow \infty$. Also, it satisfies
\begin{equation}
\label{besselDiff}
\bessel['']{\mu}{2 \pi y}=\frac{\left((2\pi y)^2+\mu^2\right)\bessel{\mu}{2\pi y}-2 \pi y \bessel[']{\mu}{2\pi y}}{(2\pi y)^2},
\end{equation}
as can be seen from the definition~\eqref{mellin-bessel} by differentiating under the integral sign and applying the well-known recursion formula for the gamma function
\begin{equation}
\label{gamma-recursion}
\Gamma(s+1) = s \Gamma(s).
\end{equation}

In certain cases, a Barnes-type integral~\eqref{barnes-integral} can be evaluated by moving the line of integration to the left or to the right and summing the residues, using the fact that the gamma function has poles at non-positive integers with residues given by
$$
\mbox{Res}(\Gamma(s), s=-n) = \frac{(-1)^n}{n!}\qquad (n = 0,1,2,\dots).
$$
When working with the resulting series, it is convenient to use the \defn{Pochhammer symbol} defined as 
\begin{equation}
\pochh{x} = x (x+1) \hdots (x+n-1)=\frac{\Gamma(x+n)}{\Gamma(x)}.
\end{equation}
Note that
\begin{equation}
\label{gammaPoch}
\pochh{1-x} = \frac{(-1)^n \Gamma(x)}{\Gamma(x-n)}.
\end{equation}

\subsection{Automorphic forms}
We recall some general facts from the theory of automorphic forms on $\gl{n}$; these facts can be found in \cite{goldfeld}.

Define the generalized upper half-plane as 
$$
\upperh{3} = \iwasawa{3}.
$$
By Iwasawa decomposition, every $z \in \upperh{3}$ can be uniquely written as $z = X Y$ with
\begin{equation}
X = \begin{pmatrix}1&x_2&x_3\\0&1&x_1\\0&0&1\end{pmatrix} \quad \text{and} \quad Y = \begin{pmatrix}y_1 y_2 &0&0\\0&y_1&0\\0&0&1\end{pmatrix}, 
\end{equation}
where $x_i,y_i \in \mathbb{R}$ and $y_i>0$. The group $\slz{3}$ acts on $\upperh{3}$ by matrix multiplication, and we are interested in functions defined on $\upperh{3}$ invariant under this action.

Let  $\nu=(\nu_1,\nu_2) \in \mathbb{C}^2$. We introduce the following three parameters (\defn{Langlands parameters})
\begin{equation}
\begin{split}
\label{sl3Parameters}
\alpha & = -\nu_1 -2 \nu_2+1,\\
\beta & = 2 \nu_1 + \nu_2-1,\\
\gamma & = -\nu_1 + \nu_2.
\end{split}
\end{equation}
Notice that $\alpha + \beta + \gamma = 0$. We will also assume that $\Re{(\alpha)} = \Re{(\beta)} = \Re{(\gamma)}=0$, that is, we assume Selberg's eigenvalue conjecture (see~\cite{goldfeld}, Conjecture 12.4.4). 

We also define
\begin{equation}
\label{lambdais}
\begin{split}
\lambda_1 & = -1-\beta \gamma - \gamma \alpha - \alpha \beta, \\
\lambda_2 & = - \alpha \beta \gamma.
\end{split}
\end{equation} 

Let $D$ be a differential operator on $\upperh{3}$ that is invariant under the action of $\slz{3}$. Any such operator can be written as a polynomial in $\Delta_1$ and $\Delta_2$, where $\Delta_1$ and $\Delta_2$ are certain differential operators given in Chapter 6 of~\cite{goldfeld} or Chapter II of~\cite{bump}. 

A \defn{Maass form of type $\nu$} for $\slz{3}$ is a smooth function $f(z)$ on $\upperh{3}$ which is
\begin{enumerate}[(a)]
\item Automorphic: $f(g z) = f(z)$ for all $g \in \slz{3}$.
\item Eigenfunction: $\Delta_i f(z)  =\lambda_i f(z)$ for $i = 1,2$ and $\lambda_i$ as in~\eqref{lambdais}.
\item Square-integrable: $\int_{\slz{3} \backslash \upperh{3}} \abs{f(z)}^2 d^*z< \infty$.
\item Cusp form: $\oneint{\oneint{\oneint{f(z)}{x_1}}{x_2}}{x_3} = 0$.
\end{enumerate}

We have the following Fourier expansion for a Maass form $f(z)$ (see~\cite{goldfeld}, Theorem 6.5.7):
\begin{equation}
\begin{split}
\label{sl3fourier}
f(z) = &  f\left(\begin{pmatrix}1&x_2&x_3\\0&1&x_1\\0&0&1\end{pmatrix}\begin{pmatrix}y_1 y_2&0&0\\0&y_1&0\\0&0&1\end{pmatrix}\right)\\
= & \sum_{(c,d)=1} \sum_{m_1=1}^\infty \sum_{m_2 \neq 0} \frac{A(m_1,m_2)}{\abs{m_1 m_2}} e^{2\pi i \left[m_1 (c x_3 + d x_1) + m_2 \Re \frac{a z_2 + b}{c z_2 + d} \right]}\\
& \qquad \qquad \times \whittaker{m_1 y_1 \abs{c z_2 +d}}{\frac{m_2 y_2}{\abs{c z_2 + d}^2}},
\end{split}
\end{equation}
where $A(m_1,m_2) \in \mathbb{C}$ are the \defn{Fourier coefficients}, $z_2=x_2 + i y_2$ and $a,b \in \mathbb{Z}$ are defined by $a d - b c = 1$, and $\whittaker{y_1}{y_2}$ is a special case of Jacquet's Whittaker function.  These higher rank Whittaker functions were introduced by Jacquet \cite{jacquet} in much greater generality; for our purposes it suffices to define \defn{Jacquet's Whittaker function} by the following explicit integral representation initially derived by Vinogradov and Takhtadzhyan \cite{takhtadzhyan} and later, in more generality, by Stade~\cite{stade}:
\begin{multline}
\label{wStade}
  \whittaker{y_1}{y_2}  =  16 (\pi y_1)^{1-\frac{\gamma}{2}} (\pi y_2)^{1+\frac{\gamma}{2}} \mellin{ \bessel{\frac{\alpha-\beta}{2}}{2 \pi y_1 \sqrt{1+u}}\\ \times \bessel{\frac{\alpha-\beta}{2}}{2 \pi y_2 \sqrt{1+u^{-1}}} }{u}{-\frac{3 \gamma}{4}}.
\end{multline}
From the decay properties of the K-Bessel function, it is clear that Jacquet's Whittaker function decays rapidly for $y_i \rightarrow \infty$, $i=1,2$. For a more precise description of its asymptotic behavior see~\cite{bump-huntley}.

Note that we suppress from notation the dependence of $\whittaker{y_1}{y_2}$ and $f(z)$ on the Langlands parameters $(\alpha, \beta, \gamma)$.  Moreover, our definition of $\whittaker{y_1}{y_2}$ differs from other definitions in the literature by a factor dependent on the values of these parameters. Since $(\alpha, \beta, \gamma)$ remain fixed throughout much of the discussion, we hope that this will not cause any confusion.  See~\cite{broughan} for a reconciliation of different definitions of Jacquet's Whittaker function as well as for the ``correct" value of the constant factor.

A Maass form $f(z)$ gives rise to an $L$-function defined by the Dirichlet series
\begin{equation}
\label{sl3Dirichlet}
L_f(s) = \sum_{n=1}^\infty \frac{A(1,n)}{n^s},
\end{equation}
which has analytic continuation to the entire complex plane and satisfies the following functional equation:
\begin{equation}
\label{sl3feq}
\Lambda_f(s) := \pi^{-\frac{3 s}{2}} \gammaf{s-\alpha}\gammaf{s-\beta} \gammaf{s-\gamma} L_f(s) = \conj{\Lambda_{f}(1-\conj{s})}.
\end{equation}
Conversely, under certain conditions (see~\cite{goldfeld}, Theorem 7.1.3) a Dirichlet series such as~\eqref{sl3Dirichlet} satisfying~\eqref{sl3feq} gives rise to a Maass form $f(z)$ for $\slz{3}$ with Langlands parameters $(\alpha,\beta,\gamma)$ and Fourier coefficients given by the following identity (\cite{bump}, Chapter 9):
$$
\sum^\infty_{m=1} \sum^\infty_{n=1} \frac{A(m,n)}{m^{s_1} n^{s_2}} = \frac{\conj{L_f(\conj{s_1})} L_f(s_2)}{\zeta(s_1+s_2)}
$$ 

The only explicitly known examples of Maass forms for $\slz{3}$ are the so-called \defn{lifts}: 
given a Maass form for $\slz{2}$, a self-dual form for $\slz{3}$ can be constructed via what's known as \defn{Gelbart-Jacquet lift} (see~\cite{goldfeld} for details). The generic forms, although more numerous, have been completely unknown until the recent numerical work by C. Bian and A. Booker \cite{booker-ams},\cite{bian} and 
D. Farmer, S. Koutsoliotas and S. Lemurell \cite{farmer-etal}, who have produced approximations for the Langlands parameters and Fourier coefficients for a number of Maass forms by computing $L$-functions and using the converse theorem above. Verifying the automorphy of these presumptive forms is one possible application of the algorithm for evaluating the Whittaker functions presented in the later section.

\subsection{Numerical integration}
To evaluate Jacquet's Whittaker function one can directly compute the defining integrals. Our method of choice for numerical integration is essentially a version of the trapezoid rule. We now recall error estimates for this method from numerical analysis.

We apply the trapezoid rule to computing the integrals of the form $\infint[-\infty]{f(x)}{x}$,
where $f(x)$ has rapid decay as $x \rightarrow \pm \infty$. The trapezoid rule in this setting gives
\begin{equation}
\label{trapezoid}
\infint[-\infty]{f(x)}{x} \approx h \sum_{k=-\infty}^{\infty} f(kh).
\end{equation}
It is a well-known fact (see, for example, \cite{fettis},~\cite{mori}, or~\cite{rubinstein}) that this approximation is extremely accurate: the discretization error of this approximation is of size $O(e^{-c/h})$. More precisely, using the Poisson summation formula, it can be shown that the error of the approximation above is given by $\sum_{k \neq 0} \hat{f}\left(\frac{k}{h}\right)$,
where $\hat{f}(y) = \infint[-\infty]{f(x) e^{-2 \pi i y x}}{x}$ denotes the Fourier transform of $f(x)$. If $\hat{f}(y)$ decays exponentially, we arrive at the error estimate given above.

Similarly, when this method is applied to the inverse Mellin transform integral
\begin{equation}
\label{poissonMult}
g(y) = \invmellin{\mellinM{g}(s) y^{-s}}{s} \approx \frac{h}{2 \pi}  \sum^\infty_{k=-\infty} \mellinM{g}(\sigma+i k h) y^{-\sigma - i k h},
\end{equation}
the error of this approximation is given by 
$$
\sum_{k\neq 0} g(y e^\frac{2 \pi k}{h}) e^{\frac{2 \pi k \sigma}{h}}. 
$$
In our applications, $g(y)$ has exponential decay as $y \rightarrow \infty$ and is bounded as $y \rightarrow 0$. One can see that the sum above decays exponentially if $k \sigma < 0$ and doubly-exponentially if $k \sigma >0$. The resulting discretization error can be made quite small by choosing an appropriate value of $\sigma$.

In our applications, the summands in the infinite sums on the right-hand side of the equations above decay exponentially. The sums can thus be easily truncated when several consecutive terms are smaller than the desired accuracy. The truncation error is then roughly of the size of the first discarded term. 

Because of the finite precision of the machine representation of real numbers, numerical results suffer from roundoff errors. We have done little to control this type of error: its effects are less pronounced when one increases working precision, and multiple precision arithmetic is readily available. We will point out which algorithms are most susceptible to roundoff error; improving those algorithms remains a research goal.

We do not work out the error estimates in more detail than given above: in practice, it suffices to know how each parameter affects the error. One then chooses parameters experimentally by comparing results obtained with one set of parameters to more accurate results obtained with a ``better'' set of parameters, or to values whose accuracy is known from theoretical considerations.

\section{Evaluating Jacquet's Whittaker function}
\label{whittaker}
The main tools for evaluating Jacquet's Whittaker functions are Stade's formula~\eqref{wStade} and the following representation in terms of a double inverse Mellin transform:
\begin{multline}
\label{wMellin}
  \whittaker{y_1}{y_2}  =  \invmmellin[(\sigma_1)(\sigma_2)]{\gammaprod{\frac{s_2-\alpha}{2},\frac{s_2-\beta}{2},\frac{s_2-\gamma}{2},\frac{s_1+\alpha}{2},\frac{s_1+\beta}{2},\frac{s_1+\gamma}{2}}{\frac{s_1+s_2}{2}} \\ \times (\pi y_1)^{1-s_1}(\pi y_2)^{1-s_2}}{s_1}{s_2}.
\end{multline}
The equivalence of the two representations for $\whittaker{y_1}{y_2}$ can be easily established by applying the double Mellin transform in $y_1$ and $y_2$ to both expressions.

Note that since $\conj{(\alpha,\beta,\gamma)} = (-\alpha,-\beta,-\gamma)$, we have
\begin{equation}
\label{whittakerDual}
\whittaker{y_2}{y_1} = \conj{\whittaker{y_1}{y_2}}.
\end{equation}

We now present several algorithms for evaluating Jacquet's Whittaker function. An implementation of all these algorithms in PARI/GP~\cite{pari} is available on the author's website \cite{mycode}. In the implementation, to avoid underflow in numerical computations, we work with $e^{\pi \abs{\alpha-\beta}}\whittaker{y_1}{y_2}$.

\subsection{Stade's formula}
\label{numStade}
We use Stade's formula~\eqn{wStade} and make a substitution $u \rightarrow e^u$ to obtain
\begin{multline*}
  \whittaker{y_1}{y_2}  =  4 (2 \pi y_1)^{1-\frac{\gamma}{2}} (2 \pi y_2)^{1+\frac{\gamma}{2}} \int_{-\infty}^{\infty} \bessel{\frac{\alpha-\beta}{2}}{2 \pi y_1 \sqrt{1+e^u}} \\ \times \bessel{\frac{\alpha-\beta}{2}}{2 \pi y_2 \sqrt{1+e^{-u}}} {e}^{-\frac{3 \gamma}{4}u} \,du.
\end{multline*}
Because of the exponential decay of the K-Bessel function, the integrand decays doubly exponentially, and the integral can be efficiently evaluated by the trapezoid rule as explained above. The only drawback of this method is the large number of evaluations of the K-Bessel function, which could be quite costly. That drawback aside, we found this method to be applicable to the widest range of values of $y_1$ and $y_2$, and used it for an informal check of results produced by other integration methods.

\subsection{Power series expansion near $(0,0)$}
\label{seriesDumb}
From equation~\eqn{wMellin}, by moving the line of integration to the left of the poles of the gamma functions and using the residue theorem (see~\cite{bump} for details), we get 
\begin{multline*}
 \whittaker{y_1}{y_2} = \sum_{(\delta_1,\delta_2,\delta_3)} (\pi y_1)^{1+\delta_1} (\pi y_2)^{1-\delta_2} \\
\times \sum_{m=0}^\infty \sum_{n=0}^\infty \gammaprod{\frac{\delta_2-\delta_1}{2}-m,\frac{\delta_2-\delta_3}{2}-m,\frac{\delta_2-\delta_1}{2}-n,\frac{\delta_3-\delta_1}{2}-n}{\frac{\delta_2-\delta_1}{2}-m-n}\\
 \times \frac{(-1)^{(m+n)} (\pi y_1)^{2 n} (\pi y_2)^{2 m}}{m! \, n!},
\end{multline*}
where the outer sum runs over the six permutations of the triple $(\alpha,\beta,\gamma)$. Using equation~\eqref{gammaPoch}, we can rewrite the expression above as follows:
\begin{multline*}
 \whittaker{y_1}{y_2} = \sum_{(\delta_1,\delta_2,\delta_3)} (\pi y_1)^{1+\delta_1} (\pi y_2)^{1-\delta_2} \gammaprod{\frac{\delta_2-\delta_3}{2},\frac{\delta_2-\delta_1}{2},\frac{\delta_3-\delta_1}{2}}{}\\
\times \sum_{m=0}^\infty \sum_{n=0}^\infty \frac{ \pochh[m+n]{1+\frac{\delta_1-\delta_2}{2}} (\pi y_1)^{2 n} (\pi y_2)^{2 m}}{\pochh[m]{1+\frac{\delta_1-\delta_2}{2}}\pochh[m]{1+\frac{\delta_3-\delta_2}{2}}\pochh[n]{1+\frac{\delta_1-\delta_2}{2}}\pochh[n]{1+\frac{\delta_3-\delta_1}{2}}  m! \, n! }.
\end{multline*}
To evaluate the coefficients of these series, one only needs six values of the gamma function. The rest of the coefficients can be easily evaluated recursively and, if desired, stored for all later calls to the Whittaker function routine.

Note that the series converges for all values of $y_1$ and $y_2$. However, for large values of $y_1$ or $y_2$ the roundoff error becomes too large for this method to be applicable. To see this, note that the summands increase when, say, $y_1$ increases, but the resulting sum decreases exponentially. Therefore, large terms must cancel out, increasing the effect of roundoff errors.

An advantage of this approach is that, given $(\alpha,\beta,\gamma)$ and $y_1$, $y_2$, one can easily determine the number of terms needed to approximate $\whittaker{y_1}{y_2}$. Once working precision is increased so that cancellation is no longer an issue, the results of the computation only suffer from truncation error, and thus can be rigorously proved to be correct to any desired accuracy.

\subsection{Power series expansion for a small argument}
\label{seriesClever}
In this section, inspired by the methods of~\cite{venkatesh}, we develop an alternative power series expansion which provides a very good way of evaluating Jacquet's Whittaker function whenever one of the arguments is small. 

From equation~\eqn{wMellin}, by moving the line of integration in $s_1$ to the left and using the residue theorem, we have
\begin{equation}
\label{threeSeries}
\begin{split}
 \whittaker{y_1}{y_2}  =&  (\pi y_1)^{1+\alpha}  \sum_{n=0}^\infty \gammaprod{\frac{\beta-\alpha}{2}-n,\frac{\gamma-\alpha}{2}-n}{}\frac{(-1)^n (\pi y_1)^{2 n}}{n!}\\
& \qquad  \qquad \qquad \times \invmellin{\gammaprod{\frac{s_2-\alpha}{2},\frac{s_2-\beta}{2},\frac{s_2-\gamma}{2}}{\frac{s_2-\alpha}{2}-n} (\pi y_2)^{1-s_2}}{s_2}\\
& + (\pi y_1)^{1+\beta} \sum_{n=0}^\infty \gammaprod{\frac{\alpha-\beta}{2}-n,\frac{\gamma-\beta}{2}-n}{}\frac{(-1)^n (\pi y_1)^{2 n}}{n!}\\
&  \qquad  \qquad \qquad \times \invmellin{\gammaprod{\frac{s_2-\alpha}{2},\frac{s_2-\beta}{2},\frac{s_2-\gamma}{2}}{\frac{s_2-\beta}{2}-n} (\pi y_2)^{1-s_2}}{s_2}\\
 & + (\pi y_1)^{1+\gamma} \sum_{n=0}^\infty \gammaprod{\frac{\alpha-\gamma}{2}-n,\frac{\beta-\gamma}{2}-n}{}\frac{(-1)^n (\pi y_1)^{2 n}}{n!}\\
&  \qquad  \qquad \qquad \times \invmellin{\gammaprod{\frac{s_2-\alpha}{2},\frac{s_2-\beta}{2},\frac{s_2-\gamma}{2}}{\frac{s_2-\gamma}{2}-n} (\pi y_2)^{1-s_2}}{s_2}.
\end{split}
\end{equation}

We focus on one of the three sums above; the same method applies to the other two after a permutation of variables. Let 
\[
I_n(y) = \invmellin{\gammaprod{\frac{s-\frac{3 \alpha}{2}}{2},\frac{s-\beta-\frac{\alpha}{2}}{2},\frac{s-\gamma-\frac{\alpha}{2}}{2}}{\frac{s-\frac{3 \alpha}{2}}{2}-n} (\pi y)^{-s}}{s}.
\]
We have, after a change of variables $s_2 \rightarrow s_2 - \frac{\alpha}{2}$ and an application of equation~\eqref{gammaPoch},
\begin{equation}
\label{seriesIn}
\begin{split}
& (\pi y_1)^{1+\alpha} \sum_{n=0}^\infty \gammaprod{\frac{\beta-\alpha}{2}-n,\frac{\gamma-\alpha}{2}-n}{}\frac{(-1)^n (\pi y_1)^{2 n}}{n!}\\
& \qquad \qquad  \qquad \qquad \times \invmellin{\gammaprod{\frac{s_2-\alpha}{2},\frac{s_2-\beta}{2},\frac{s_2-\gamma}{2}}{\frac{s_2-\alpha}{2}-n} (\pi y_2)^{1-s_2}}{s_2}\\
& =  (\pi y_1)^{1+\alpha} \gammaprod{\frac{\beta-\alpha}{2},\frac{\gamma-\alpha}{2}}{} (\pi y_2)^{1+\frac{\alpha}{2}} \sum_{n=0}^\infty \frac{(-1)^n (\pi y_1)^{2 n}  I_n(y_2)}{\pochh{1+\frac{\alpha-\beta}{2}} \pochh{1+\frac{\alpha-\gamma}{2}} n!}.
\end{split}
\end{equation}
The integrals $I_n(y)$ satisfy a recurrence relation,  
\begin{equation}
\label{InRecurse}
\begin{split}
I_{n+1}(y) & = \invmellin{\gammaprod{\frac{s-\frac{3 \alpha}{2}}{2},\frac{s-\beta-\frac{\alpha}{2}}{2},\frac{s-\gamma-\frac{\alpha}{2}}{2}}{\frac{s-\frac{3 \alpha}{2}}{2}-n-1} (\pi y)^{-s}}{s}\\
& = \invmellin{\gammaprod{\frac{s-\frac{3 \alpha}{2}}{2},\frac{s-\beta-\frac{\alpha}{2}}{2},\frac{s-\gamma-\frac{\alpha}{2}}{2}}{\frac{s-\frac{3 \alpha}{2}}{2}-n} \left(\frac{s-\frac{3\alpha}{2}}{2}-n-1\right) (\pi y)^{-s}}{s}\\
& = -\frac{1}{2} \left(y I_n'(y) + \left(\frac{3\alpha}{2}+2n+2\right) I_n(y) \right).
\end{split}
\end{equation}
In addition,
\begin{equation}
\label{I0}
I_0(y)  = \invmellin{\gammaprod{\frac{s-\beta-\frac{\alpha}{2}}{2},\frac{s-\gamma-\frac{\alpha}{2}}{2}}{} (\pi y)^{-s}}{s} = 4 \bessel{\frac{\beta-\gamma}{2}}{2 \pi y},
\end{equation}
where we used~\eqref{mellin-bessel} and the fact that $\alpha+\beta+\gamma=0$.

Using~\eqref{I0}, the recursion for $I_n$ given in~\eqref{InRecurse}, 
 and the differential equation~\eqref{besselDiff} for the K-Bessel function, we can write
\[
I_n(y) = (-2)^{-n} \left(P_n(y) \bessel{\mu}{2 \pi y}+ 2 \pi y Q_n(y) \bessel[']{\mu}{2 \pi y}\right),
\]
with $\mu = \frac{\beta-\gamma}{2}$, and polynomials $P_n(y)$ and $Q_n(y)$ defined recursively by 
\begin{align*}
P_{n+1}(y) & = y P_n'(y)+ \left((2\pi y)^2+\mu^2\right) Q_n(y)+a_n P_n(y) & P_0(y)&=4\\
Q_{n+1}(y) &= P_n(y) + y Q_n'(y) + a_n Q_n(y) & Q_0(y) & =0,
\end{align*}
where $a_n=\frac{3\alpha}{2}+2n+2.$

Substituting this result in~\eqref{seriesIn}, we get 
\begin{multline*}
(\pi y_1)^{1+\alpha} \gammaprod{\frac{\beta-\alpha}{2},\frac{\gamma-\alpha}{2}}{} (\pi y_2)^{1+\frac{\alpha}{2}} \sum_{n=0}^\infty \frac{(-1)^n (\pi y_1)^{2 n}}{\pochh{1+\frac{\alpha-\beta}{2}} \pochh{1+\frac{\alpha-\gamma}{2}} n!} I_n(y_2)\\
 = (\pi y_1)^{1+\alpha} \gammaprod{\frac{\beta-\alpha}{2},\frac{\gamma-\alpha}{2}}{} (\pi y_2)^{1+\frac{\alpha}{2}}\\ \times \sum_{n=0}^\infty \frac{\left(P_n(y_2) \bessel{\mu}{2 \pi y_2}+ 2 \pi y Q_n(y_2) \bessel[']{\mu}{2 \pi y_2}\right)(\pi y_1)^{2 n}}{ \pochh{1+\frac{\alpha-\beta}{2}} \pochh{1+\frac{\alpha-\gamma}{2}} 2^n n!},
\end{multline*}
with $P_n(y)$, $Q_n(y)$ as above. Similar expressions are obtained for the other two sums in~\eqref{threeSeries}.

These polynomials $P$ and $Q$ can be computed once and stored for each subsequent evaluation of the function. Then the Whittaker function can be quickly and accurately evaluated with six calls to a K-Bessel routine and evaluation of polynomials.

This method is very effective for small values of $y_2$; if instead $y_1 < y_2$ is small, one should apply equation~\eqref{whittakerDual} before using this algorithm. In fact, whenever either one of the arguments of $\whittaker{y_1}{y_2}$ is of a reasonable size, we found this method to be the best way of evaluating Jacquet's Whittaker function. It avoids numerical integration and uses power series in only one variable. The problems with this method start to surface when both arguments are large: we run into the cancellation problems just like in Section~\ref{seriesDumb}. For large arguments we resort to Stade's formula from Section~\ref{numStade}.

\subsection{Numerical inverse Mellin transform}
\label{numMellin}
When one evaluates a Maass form for $\slz{3}$ using the Fourier expansion~\eqref{sl3fourier}, the Whittaker function has arguments of a very specific form: one repeatedly computes $\whittaker{y_1 \abs{z}}{\frac{y_2}{\abs{z}^2}}$ for fixed $y_1$, $y_2$ and varying $z$. Now let $D = y_1^2 y_2$ and define $\whittaker[*]{D}{y_2} = \whittaker{y_1}{y_2}$. Then 
$$
\whittaker{y_1 \abs{z}}{\frac{y_2}{\abs{z}^2}} = \whittaker[*]{D}{\frac{y_2}{\abs{z}^2}}.
$$
When $z$ changes, $D$ stays fixed, effectively making $W^*$ a single-variable function. We seek to exploit this fact in our implementation.

From equation~\eqn{wMellin}, we have 
\begin{equation*}
\begin{split}
  \whittaker{y_1}{y_2}  & = \frac{1}{4 \pi^2} \iinfint{\gammaprod{\frac{\sigma_2 + i t_2 -\alpha}{2},\frac{\sigma_2 + i t_2-\beta}{2},\frac{\sigma_2 + i t_2-\gamma}{2},\frac{\sigma_1 + i t_1+\alpha}{2},\frac{\sigma_1 + i t_1+\beta}{2},\frac{\sigma_1 + i t_1+\gamma}{2}}{\frac{\sigma_1 + i t_1+\sigma_2 + i t_2}{2}}\\
& \qquad \times (\pi y_1)^{1-\sigma_1- i t_1} (\pi y_2)^{1-\sigma_2- i t_2} }{t_1}{t_2}\\
& = \frac{(\pi^3 D)^\frac{1-\sigma_1}{2} (\pi y_2)^\frac{1-2 \sigma_2 + \sigma_1}{2}}{4 \pi^2} \\
& \qquad \times \iinfint{%
\gammaprod{\frac{\sigma_2 + i t_2 -\alpha}{2},\frac{\sigma_2 + i t_2-\beta}{2},\frac{\sigma_2 + i t_2-\gamma}{2}}{}\\
&\qquad  \times %
\gammaprod{\frac{\sigma_1 + i t_1+\alpha}{2},\frac{\sigma_1 + i t_1+\beta}{2},\frac{\sigma_1 + i t_1+\gamma}{2}}{\frac{\sigma_1 + i t_1+\sigma_2 + i t_2}{2}}%
  (\pi^3 D)^{\frac{- i t_1}{2}} (\pi y_2)^{- i t_2+\frac{i t_1}{2}} }{t_1}{t_2}.
 \end{split}
\end{equation*}
After a change of variables, we obtain
\begin{equation}
\begin{split}
\label{whittakerD}
&\whittaker{y_1}{y_2} = \whittaker[*]{D}{y_2} \\
& \qquad = \frac{(\pi^3 D)^\frac{1-\sigma_1}{2} (\pi y_2)^\frac{1-2 \sigma_2 + \sigma_1}{2}}{2 \pi^2} \\
& \qquad \qquad \times \iinfint{\gammaprod{\frac{\sigma_1+\alpha}{2}+ i t_1,\frac{\sigma_1 +\beta}{2}+ i t_1,\frac{\sigma_1 +\gamma}{2}+ i t_1}{\frac{\sigma_1 +\sigma_2 + i (t_2+ 3 t_1)}{2}} \\
& \qquad \qquad \times \gammaprod{\frac{\sigma_2 + i (t_2+t_1) -\alpha}{2},\frac{\sigma_2 + i (t_2+t_1)-\beta}{2},\frac{\sigma_2 + i (t_2+t_1)-\gamma}{2}}{}\\
& \qquad \qquad \times  (\pi^3 D)^{- i t_1} (\pi y_2)^{- i t_2}}{t_1}{t_2}.
\end{split}
\end{equation}

We have written $\whittaker{y_1}{y_2}$ as a function of $D = y_1^2 y_2$ and $y_2$, and now exploit this representation in the implementation as follows. Given the parameters $(\alpha,\beta,\gamma)$, one pre-computes all the gamma factors needed to evaluate the integral in~\eqref{whittakerD} by the trapezoid rule:
\begin{equation*}
\begin{split}
\whittaker{y_1}{y_2} & \approx \frac{(\pi^3 D)^\frac{1-\sigma_1}{2} (\pi y_2)^\frac{1-2 \sigma_2 + \sigma_1}{2}}{2 \pi^2} h_1 h_2\\
& \qquad \times \sum_{k_2=-N_2}^{N_2} \left\{\sum_{k_1=-N_1}^{N_1} \gammaprod{\frac{\sigma_1+\alpha}{2}+ i k_1 h_1,\frac{\sigma_1 +\beta}{2}+ i k_1 h_1,\frac{\sigma_1 +\gamma}{2}+ i k_1 h_1}{\frac{\sigma_1 +\sigma_2 + i (t_2+ 3 k_1 h_1)}{2}} \right.\\
& \qquad \qquad \times \gammaprod{\frac{\sigma_2 + i (k_2 h_2+k_1 h_1) -\alpha}{2},\frac{\sigma_2 + i (k_2 h_2+k_1 h_1)-\beta}{2}}{}\\
& \qquad \qquad \times \left. \gammaprod{\frac{\sigma_2 + i (k_2 h_2+k_1 h_1)-\gamma}{2}}{}(\pi^3 D)^{- i k_1 h_1} \right\} (\pi y_2)^{- i k_2 h_2}.
\end{split}
\end{equation*}
For some choice of parameters $h_1$, $h_2$, $\sigma_1$, $\sigma_2$, $N_1$ and $N_2$ this will result in a very good approximation, see~\eqref{poissonMult} and the discussion that follows.

Now, given a value of $D$, the inner sum is evaluated for each value of $k_2$ and the results are stored. For each subsequent evaluation of Jacquet's Whittaker function with the same value of $D$ and varying $y_2$, one only needs to re-evaluate the outer sum --- a very fast operation.

The main drawback of this method is that it is prone to roundoff errors, especially for larger values of $y_2$. This can be seen quite easily: the Whittaker function is exponentially decaying in $y_2$, but the integrand can be quite large, especially for larger values of $\abs{\alpha}$, $\abs{\beta}$ and $\abs{\gamma}$. Thus enormous cancellation must occur, leading to inevitable loss of precision. Still, we have found this approach to be the fastest of all the ones we explored. We increase working precision and use this method to evaluate Maass forms for $\slz{3}$ in the next section. 


\section{Evaluating Maass forms for $\slz{3}$}
\label{maass}
Suppose we are given the Langlands parameters $(\alpha,\beta,\gamma)$ and the Fourier coefficients $A(m_1,m_2)$ of a Maass form $f(z)$ for $\slz{3}$. The problem is to evaluate it at some point $z_0 \in \upperh{3}$.

Since every Maass form for $\slz{3}$ is even (see~\cite{goldfeld}), we can rewrite equation~\eqref{sl3fourier} as
\begin{equation}
\begin{split}
\label{sl3cosine}
f(z) & = f\left(\begin{pmatrix}1&x_2&x_3\\0&1&x_1\\0&0&1\end{pmatrix}\begin{pmatrix}y_1 y_2&0&0\\0&y_1&0\\0&0&1\end{pmatrix}\right)\\
& = 4 \sum_{m_1=1}^\infty \sum_{m_2=1}^\infty \frac{A(m_1,m_2)}{m_1 m_2} \Biggl[ \cospi{m_2 x_2}\cospi{m_1 x_1} \whittaker{m_1 y_1}{m_2 y_2} \\%
& \qquad \qquad + \sum_{\substack{c \geq 1\\(c,d)=1}}\cospi{m_1 (c x_3 + d x_1)} \cospi{\frac{m_2}{c} \left(a-\Re{\frac{1}{c z_2 + d}}\right)}\\
& \qquad \qquad \qquad \qquad \qquad \times \whittaker{m_1 y_1 \abs{c z_2+d}}{\frac{m_2 y_2}{\abs{c z_2+d}^2}}\Biggr],
\end{split}
\end{equation}
where $z_2=x_2 + i y_2$, $a \in \mathbb{Z}$ is defined by $a d\equiv 1\pmod c$ and $\whittaker{y_1}{y_2}$ denotes Jacquet's Whittaker function.

Now, given $(\alpha,\beta,\gamma)$ and an accuracy goal $\epsilon$, one can determine $C = C(\alpha,\beta,\gamma)$ such that $\abs{\whittaker{y_1}{y_2}} < \epsilon$ whenever $y_1>C$ or $y_2>C$. For example, $C$ can be easily determined numerically by evaluating the Whittaker function for increasing values of $y_1$. Once $C$ is known, one can truncate the inner sum in~\eqref{sl3cosine}: for a fixed pair $c$, $d$, a term contributes whenever 
$$m_1 y_1 \abs{c z_2+d} < C \quad \text{ and } \quad \frac{m_2 y_2}{\abs{c z_2 +d}^2} <C.$$
Therefore, $c$ and $d$ must satisfy
\begin{equation}
\label{cdcutoff}
\sqrt{\frac{m_2 y_2}{C}}< \abs{c z_2 + d} < \frac{C}{m_1 y_1}.
\end{equation}
It follows that for each $m_1 y_1$, $m_2 y_2$ and $z_2$, the infinite sum over all relatively prime integers $c$ and $d$ with $c\geq 1$ can be truncated whenever the condition~\eqref{cdcutoff} is not satisfied. Because of the exponential decay of Jacquet's Whittaker function, truncation of the outer sums over $m_1$ and $m_2$ does not pose a problem.

To evaluate the inner sum for fixed values of $m_1 y_1$, $m_2 y_2$, and $z_2$, we use the algorithm from Section~\ref{numMellin}: $D = (m_1 y_1)^2 m_2 y_2$ is fixed, and we compute the value of $\whittaker{m_1 y_1 \abs{c z_2+d}}{\frac{m_2 y_2}{\abs{c z_2+d}^2}}$ for all $c$, $d$ that satisfy~\eqref{cdcutoff}. As mentioned in Section~\ref{numMellin}, when $D$ is fixed, each Whittaker function evaluation takes very little time. The only time-consuming operation is changing the value of $D$, but that happens much less frequently during the computation. 

We apply this algorithm to evaluate several $\slz{3}$-Maass forms $f(z)$ at $z_0 = I_3$, the $3\times3$ identity matrix. The Maass forms are lifts of Maass forms for $\slz{2}$, which have been computed to high precision in~\cite{venkatesh}. The Langlands parameters of lifted forms are $(-2 i r, 2 i r, 0)$, where $r$ is the type of the $\slz{2}$-Maass form as in~\cite{venkatesh}. The time required for each computation and the number of the Fourier coefficients needed to evaluate $f(z_0)$ with error not exceeding $10^{-12}$ are summarized in the table below.
\begin{table}[ht]
\centering
\begin{tabular}{|c|c|c|}\hline
\label{table:maass}
 $r$ & Coeffs & Time (s) \\ \hline
 $9.533695$ & 120 & $94$ \\ \hline
 $13.779751$ & 203 & $204$ \\ \hline
 $17.738563$ & 379 & $440$ \\ \hline
 $22.194674$ & 590 & $972$ \\ \hline
 $26.056918$ & 861 & --- \\ \hline
 $35.431665$ & 1806 & --- \\ \hline
\end{tabular}
\caption{Time and number of the Fourier coefficients needed to evaluate a lifted Maass form on a Pentium D 2.80GHz.}
\end{table}

The number of coefficients is essentially the maximum value of $m_2$ for which we get a non-negligible contribution from the innermost sum in~\eqref{sl3cosine}. As should be expected from~\eqref{cdcutoff}, the maximum value of $m_1$ resulting in non-negligible contribution is much smaller, not exceeding $10$ for the first cases. It would be interesting to see if there is another way to evaluate a Maass form that requires fewer coefficients for the same level of precision, perhaps by making the expansion more symmetric in $m_1$ and $m_2$. Given that Farmer et al \cite{farmer-etal} require only about $20$ coefficients to locate a candidate Maass form, it is somewhat surprising that evaluating such a form would require an order of magnitude more coefficients.
 
In another application of this algorithm, we evaluate a generic Maass form with Langlands parameters $(-14.141638, -2.380388, 16.522027)$ discovered by Ce Bian \cite{bian}.  Let
\begin{eqnarray*}
S_1 = \begin{pmatrix} 1 & 0 & 0 \\ 0 & 0 & -1 \\ 0 & 1 & 0 \end{pmatrix} & \text{and} & S_2 = \begin{pmatrix} 0 & -1 & 0 \\ 1 & 0 & 0 \\ 0 & 0 & 1 \end{pmatrix}
\end{eqnarray*}
be two elements of $\slz{3}$. We pick several points $z_0=((x_1,x_2,x_3),(y_1,y_2)) \in \upperh{3}$ and evaluate $f(z_0)$ and $f(\gamma z_0)$ where $\gamma$ is a word in $S_i$'s. Some results are given in the table below.

\begin{table}[ht]
\centering
\begin{tabular}{|c|c|}\hline
\label{table:whittaker}
$z_0$ & $f(z_0) = f(x_1,x_2,x_3,y_1,y_2)$ \\ \hline
$((0,0,0),(0.9,0.9))$ & $-79.779900 - 0.000044759125i$\\
$S_1 S_2 S_1 ((0,0,0),(0.9,0.9))$ & $-79.780569 - 0.000025792266i$\\ \hline
$((0,0,0),(0.9,10/9))$ & $-168.098614 + 7.88405197i$\\
$S_2 S_1 ((0,0,0),(0.9,10/9))$ & $-168.099135 + 7.88392995i$\\ \hline
$((0.2,0.4,-0.1),(0.95,0.9))$ & $-16.1067874 + 12.9212510i$\\
$S_2 S_1 ((0.2,0.4,-0.1),(0.95,0.9))$ & $-16.1066031 + 12.9212690i$\\ \hline
$((0.726,0.325,0.983),(0.541,0.578))$ & $11.3789573 + 0.0986941489i$\\
$S_2 S_1 ((0.726,0.325,0.983),(0.541,0.578))$ & $11.3788905 + 0.0986972163i$\\ \hline
\end{tabular}
\caption{Values of a Maass form $f(z)$ with $(\alpha,\beta,\gamma) = (-14.141638, -2.380388, 16.522027)$ discovered by Bian and Booker}
\end{table}

The agreement between the values of $f$ at two $\slz{3}$-equivalent points of $\upperh{3}$ is quite good. Given the precision of the Langlands parameters and the Fourier coefficients, the errors of this size are to be expected. We thus have further evidence of automorphy of Bian's first example. A rigorous numerical verification of automorphy along the lines of \cite{venkatesh} remains a research goal. 

\section*{Acknowledgments}
The advice of Professor Dorian Goldfeld is gratefully acknowledged. I also thank Alex Kontorovich, David Farmer, Andrew Booker and other participants of AIM workshop on computing arithmetic spectra for helpful discussions. Special thanks are due to Ce Bian for providing the data for the numerical experiments.

\addcontentsline{toc}{chapter}{Bibliography}
\def\cprime{$'$} \def\cprime{$'$}

\end{document}